\documentclass{amsart}
\begin{document}
\title{The Harmonic Series and the $n$th Term Test for Divergence}
\author{David M.~Bradley}
\address{Department of Mathematics and Statistics\\
 University of Maine \\
        Orono, Maine 04469-5752\\ U.S.A.}
\email{bradley@math.umaine.edu; dbradley@member.ams.org}
\date{February 29, 2000}
\maketitle

The harmonic series $\sum_{n=1}^\infty 1/n$ is a popular example
of a divergent series whose terms tend to zero.  Another example
is $\sum_{n=1}^\infty \log(1+1/n),$ whose   partial sums are
unbounded because they telescope to $\log n$:
\[
   \sum_{k=1}^{n-1}\log\bigg(1+\frac{1}{k}\bigg)=
   \sum_{k=1}^{n-1}\log\bigg(\frac{k+1}{k}\bigg)=
   \sum_{k=1}^{n-1}\big(\log(k+1)-\log k\big)=
   \log n.
\]

This telescoping series can be used to show that the harmonic
series diverges. Start with the classical inequality $x\ge
\log(1+x)$, which is valid for all $x>-1$ (compare the derivatives
of both sides and integrate from $0$ to $x$).  Now put $x=1/k$
with $k=1,2,3,\ldots,n-1$ and add the resulting inequalities,
obtaining
\[
   \sum_{k=1}^{n-1}\frac{1}{k}
   \ge \sum_{k=1}^{n-1}\log\bigg(1+\frac{1}{k}\bigg)
   = \log n.
\]
\end{document}